\documentclass[10pt]{article}
\usepackage{color}
\usepackage{amssymb, amsmath}
\usepackage{epsfig}
\usepackage{graphicx}

\newtheorem{theorem}{Theorem}

\newtheorem{lemma}{Lemma}%

\newcommand{\proof}{\noindent\textbf{Proof.~}}

\newcommand{\qed}{\space\hfill\hspace*{\fill} $\vbox{\hrule\hbox{\vrule
height1.3ex\hskip1.3ex\vrule}\hrule}$\hss\vskip\topsep\relax}

\begin{document}

\title{A novel approach to construct numerical methods for stochastic differential equations}

\author{Nikolaos Halidias \\
{\small\textsl{Department of Statistics and Actuarial-Financial Mathematics }}\\
{\small\textsl{University of the Aegean }}\\
{\small\textsl{Karlovassi  83200  Samos, Greece} }\\
{\small\textsl{email: nikoshalidias@hotmail.com}}}

\maketitle

\begin{abstract}In this paper we propose a new numerical method
for solving stochastic differential equations (SDEs).  As an
application of this method we propose an explicit numerical scheme
for a super linear SDE for which the usual Euler scheme diverges.
\end{abstract}

{\bf Keywords:}   Explicit numerical scheme, super linear
stochastic differential equations.

{\bf AMS subject classification:}  60H10, 60H35.

\section{Introduction}
  Let $(\Omega, {\cal F},
\mathbb{P}, {\cal F}_t)$ be a complete probability space with a
filtration and let  a Wiener process $(W_t)_{t \geq 0}$ defined on
this space. Consider the following stochastic differential
equation,
\begin{eqnarray}
x_t = x_0 + \int_0^t a(x_s)ds + \int_0^t b(x_s)dW_s,
\end{eqnarray}
where $a,b: \mathbb{R}_+ \times \mathbb{R} \to \mathbb{R}$ are
measurable functions
 and $x_0$ such that is ${\cal F}_0$-measurable and square
 integrable.
Suppose that this problem has a unique strong solution $x_t$.

{\em {\bf Assumption A}  Assume that there exist $f(x,y)$, $g(x,y):
\mathbb{R} \times \mathbb{R} \to \mathbb{R}$ such that $f(x,x) =
a(x)$, $g(x,x) = b(x)$. Furthermore, assume that for any $R > 0$ and
any
 $|x_1|,|x_2|,|y_1|,|y_2| \leq R$  we have
\begin{eqnarray*}
|f(x_1,y_1) - f(x_2,y_2)| + |g(x_1,y_1) -
g(x_2,y_2)|  \leq C_R (|x_1-x_2| + |y_1-y_2|),
\end{eqnarray*}}
for some $C_R$ depending on $R$ and on $a,b,f,g$.

Let $0 = t_0 < t_1 < ...<t_n = T$ and set $\Delta = \frac{T}{n}$.
 For any $t \in [t_{k},t_{k+1}]$ consider the following stochastic differential equation,
\begin{eqnarray}
y_t = y_{t_{k}}+ \int_{t_k}^{t}
f(y_s,y_{t_k})ds + \int_{t_k}^{t}
g(y_s,y_{t_k}) dW_s, \quad t \in [t_k,t_{k+1}]
\end{eqnarray}
where we assume that  SDE  (2), in each step, has a unique strong
solution which we denote by $y_t$ and $y_0 = x_0$. Note, that the
coefficients in (2) are random. Let us call (2) a ''semi-discrete''
numerical scheme because, in general, we   discretize  only a  part
of our original SDE. In practice we will choose $f,g$ such that this
numerical scheme will have a known explicit solution, for example,
the ''semi-discrete''  SDE (2) may be linear. We can write the above
numerical scheme more compactly,
\begin{eqnarray}
y_t = x_0 + \int_0^t f(y_s,y_{\hat{s}})ds + \int_0^t
g(y_s,y_{\hat{s}})dW_s,
\end{eqnarray}
where
\begin{eqnarray*}
\hat{s} = t_k, \;\;\;   \quad s \in [t_k,t_{k+1}].
\end{eqnarray*}
The solution $y_t$ of (3) depends on $\Delta$ and we should use a
notation like this, $y_t^{\Delta}$, but we don't do here  for
simplicity. Note that we use the  second variable in $f,g$ to denote
the discretized part of the original SDE.

The SDE (3), i.e. the ''semi-discrete'' SDE, is not an implicit
numerical scheme, because in order to solve for $y_t$ we have to
solve a stochastic differential equation and not an algebraic
equation. In our setting, we can not  reproduce the known implicit
numerical schemes but we can reproduce for example the  Euler
scheme. So the usual Euler scheme belongs to our setting choosing
$f(x,y) = a(y)$ and $g(x,y) = b(y)$. Another interesting way to
choose $f,g$ is $f(x,y) = -\frac{1}{2} b^{'}(y)b(y) + a(y)
+\frac{1}{2} b^{'}(x)b(x) $ and $g(x,y) = b(x)$ (see
\cite{Halidias2}). This comes from the fact that the following SDE,
\begin{eqnarray*}
x_t = x_0 + \int_0^t \frac{1}{2} b^{'}(x_s) b(x_s)ds + \int_0^t
b(x_s)dW_s,
\end{eqnarray*}
has a known explicit solution (see \cite{Kloeden}, p. 117). We can
use also other, more sophisticated, SDEs with known explicit
solutions as these described in \cite{Kloeden}. Thus, our method
here is more general than \cite{Halidias2} because we can arrive to
other SDEs, like linear SDES, choosing suitable $f,g$. The main
advantage of our method is that we produce always explicit numerical
schemes in contrast to other interesting but implicit methods (see
for example \cite{Mao4}, \cite{Mao3}).

Another suitable choice, that we are going to use in our example in
this paper, is $f(x,y) = \frac{a(y)}{y} x$ and $g(x,y) =
\frac{b(y)}{y}x$. Then the resulting ''semi-discrete'' SDE will be a
linear stochastic differential  equation with known explicit
solution.

Let us point out that our main result and setting remains true for
the multidimensional case but in order to present our method as
simply as possible we will remain in the one dimensional case.

 In the second section we will
state and prove a convergence result and next we will give an
application. For a general study of the numerical analysis of
stochastic differential equations one can see \cite{Kloeden}.

\section{Convergence of the semi discrete numerical scheme}

In this section we shall  prove  that our numerical scheme
converges to the true solution.
\begin{theorem}
Assume that Assumption A holds and that (2) has a unique strong solution in
each interval. Suppose that
\begin{eqnarray*}
\mathbb{E} ( \sup_{0 \leq t \leq T} |x_t|^{p})  < A, \quad
\mathbb{E} ( \sup_{0 \leq t \leq T} |y_t|^{p}) < A
\end{eqnarray*}
for some $p > 2$, $A>0$ independent of $\Delta$. Then
\begin{eqnarray*}
\lim_{\Delta \to 0} \mathbb{E} \sup_{0 \leq t \leq T} |y(t)
- x(t)|^2 = 0.
\end{eqnarray*}
\end{theorem}

\proof
 Set $ \rho_R = \inf \{ t \in [0,T]:
 |x_t| \geq R \}$ and $ \tau_R = \inf \{ t \in [0,T] : |y_t| \geq R
\}$. Let $\theta_R = \min \{ \tau_R, \rho_R \}$.

At first, let us estimate the following probability
\begin{eqnarray*}
\mathbb{P}(\tau_R \leq T) = \mathbb{E} \left[ \mathbb{I}_{\{
\tau_R \leq T \} } \frac{|y_{\tau_R}^p|}{R^p} \right]  \leq
\frac{A}{R^p}.
\end{eqnarray*}

Therefore we can prove that $\mathbb{P}(\tau_r \leq T \mbox{ or }
\rho_R \leq T ) \leq \frac{2 A}{R^p}$. Using Young inequality we
obtain, for any $\delta > 0$,
\begin{eqnarray*}
\mathbb{E}\left(\sup_{0 \leq t \leq T } |y_t -x_t|^2\right) \leq
\mathbb{E}\left( \sup_{ 0 \leq t \leq T} \left|y_{t
\wedge\theta_R} - x_{t \wedge\theta_R}\right|^2\right) +
\frac{2^{p+1} \delta A}{p} + \frac{(p-2)2A}{p
\delta^{\frac{2}{p-2}}R^p}.
\end{eqnarray*}

We shall estimate now the term $|x_{t\wedge \theta_R} - y_{t\wedge
\theta_R}|^2$ as follows, using  Assumption A for
$f(\cdot,\cdot)$,
\begin{eqnarray*}
& & |x_{t\wedge \theta_R} - y_{t\wedge \theta_R}|^2 = \\
 & & \left| \int_0^{t\wedge \theta_R}
(f(x_s,x_s)-f(y_s,y_{\hat{s}}))  ds +
\int_0^{t\wedge \theta_R}  (g(x_s,x_s)-g(y_s,y_{\hat{s}}))dW_s \right|^2 \\
& & \leq 2 \int_0^{t\wedge \theta_R}
|f(x_s,x_s)-f(y_s,y_{\hat{s}})|^2   ds + 2 |\int_0^{t\wedge
\theta_R} (g(x_s,x_s)-g(y_s,y_{\hat{s}}))dW_s|^2 \\
& & \leq C_R \int_0^{t\wedge \theta_R} (|x_s-y_s|^2 +
|x_s-y_{\hat{s}}|^2)ds + 2|\int_0^{t\wedge \theta_R}
(g(x_s,x_s)-g(y_s,y_{\hat{s}}))dW_s|^2
\end{eqnarray*}
Note that $C_R$ will be different from line to line. We can write
now,
\begin{eqnarray*}
\sup_{0 \leq t \leq s} |x_{t\wedge \theta_R} - y_{t\wedge
\theta_R}|^2 \leq C_R \int_0^s (|x_{r\wedge \theta_R}-y_{r\wedge
\theta_R}|^2 + |x_{r\wedge \theta_R} - y_{\widehat{r\wedge
\theta_R}}|^2)dr +\\ 2 \sup_{0 \leq t \leq s} |\int_0^{t\wedge
\theta_R} (g(x_s,x_s)-g(y_s,y_{\hat{s}}))dW_s|^2.
\end{eqnarray*}
Taking expectations on both sides, using  Doob's martingale
inequality for the second term at the right hand side and
Assumption A for $g(\cdot,\cdot)$ we arrive at,
\begin{eqnarray*}
\mathbb{E} (\sup_{0 \leq t \leq s} |x_{t\wedge \theta_R} -
y_{t\wedge \theta_R}|^2) \leq C_R \mathbb{E} \int_0^s (|x_{r\wedge
\theta_R}-y_{r\wedge \theta_R}|^2 + |x_{r\wedge \theta_R} -
y_{\widehat{r\wedge \theta_R}}|^2)dr \leq \\ C_R \int_0^s
\mathbb{E} \sup_{ 0 \leq l \leq r} |x_{l\wedge
\theta_R}-y_{l\wedge \theta_R}|^2 dr + C_R \int_0^s \mathbb{E}
|y_{r\wedge \theta_R} - y_{\widehat{r\wedge \theta_R}}|^2 dr.
\end{eqnarray*}

We shall estimate the term $\mathbb{E}|y_{t\wedge \theta_R} -
y_{\widehat{{t\wedge \theta_R}}}|^2$. We begin with,
\begin{eqnarray*}
|y_{t\wedge \theta} - y_{\widehat{t\wedge \theta_R}}|^2 = |
\int_{{\widehat{t\wedge \theta_R}}}^{t\wedge \theta_R}
f(y_s,y_{\hat{s}})ds + \int_{{\widehat{t\wedge
\theta_R}}}^{t\wedge \theta_R} g( y_s,y_{\hat{s}})dW_s |^2 \leq
\\ 2 (\int_{{\widehat{t\wedge \theta_R}}}^{t\wedge \theta_R} f(y_s,y_{\hat{s}})ds)^2 +
2 |\int_{{\widehat{t\wedge \theta_R}}}^{t\wedge \theta_R}
g(y_s,y_{\hat{s}})dW_s|^2.
\end{eqnarray*}
Taking expectations, using Ito's isometry and  the fact that
$|f(y_s,y_{\hat{s}})|, |
g(y_s,y_{\hat{s}})| \leq C_R$ we have that,
\begin{eqnarray*}
\mathbb{E} |y_{t\wedge \theta} - y_{\widehat{t\wedge \theta}}|^2
\leq C_R \Delta.
\end{eqnarray*}

Collecting all together and  using  the Gronwall inequality we
arrive at,
\begin{eqnarray*}
\mathbb{E}\left(\sup_{0 \leq t \leq T } |y_t - x_t|^2\right) \leq
C_R \Delta + \frac{2^{p+1} \delta A}{p} + \frac{(p-2)2A}{p
\delta^{\frac{2}{p-2}}R^p}.
\end{eqnarray*}

Now, given $\varepsilon > 0$ choose $\delta > 0$ such that
$\frac{2^{p+1} \delta A}{p} < \varepsilon/3$ and then choose $R$
such that $\frac{(p-2)2A}{p \delta^{\frac{2}{p-2}}R^p} <
\varepsilon/3$. Finally, choose $\Delta$ small enough to get the
desired result.

\qed

\section{Example}

Below we will propose an explicit numerical scheme for a super
linear SDE as an application of Theorem 1. Consider the following
SDE,
\begin{eqnarray*}
x_t = x_0 - \int_0^t x^3_sds + \int_0^t b x_s dW_s,
\end{eqnarray*}
with $x_0 \in \mathbb{R}_+$ and $b \in \mathbb{R}$. This SDE has a
unique strong solution and we know that the usual Euler scheme
diverges (see \cite{Kloeden1}).

We propose the following ''semi-discrete'' numerical scheme which
is explicit,
\begin{eqnarray*}
y_t = x_0 + \int_0^t y_s (-y_{\hat{s}}^2) ds + \int_0^t b y_s dW_s.
\end{eqnarray*}
This semi discrete scheme is a linear SDE with unique strong
solution,
\begin{eqnarray*}
y_t =  x_0 e^{-\int_0^t (y_{\hat{s}}^2 + \frac{b^2}{2})ds +b W_t}.
\end{eqnarray*}

We need the moment bounds which we will prove in the next lemma.
\begin{lemma}
Suppose that  $x_0 > 0$ and $x_0 \in \mathbb{R}$. Then there
exists some $A
> 0$ such that
\begin{eqnarray*}
\mathbb{E}  ( \sup_{0 \leq t \leq T} |y_t|^p) < A,  \quad
\mathbb{E}( \sup_{0 \leq t \leq T} |x_t|^p) < A, \quad \mathbb{E}(
\sup_{0 \leq t \leq T} \frac{1}{|x_t|^2}) < A,
\end{eqnarray*}
for any $p \geq 2$.
\end{lemma}

\proof  Set $r= \min \{ \rho, \tau \}$ where $\rho = \inf \{ t \in
[0,T] : |x_t| > R \}$ and $\tau = \inf \{ t \in [0,T] : |y_t|
> R \}$ for some $R > 0$.

 Using Ito's formula on $|y_{t \wedge r}|^p$ we obtain,
\begin{eqnarray*}
 |y_{t \wedge r}|^ p = |x_0|^p + \int_0^t (p  |y_s|^p
(-y_{\hat{s}}^2) +\frac{b^2 p (p-1)}{2} |y_s|^{p}) \mathbb{I}_{ \{ (0,r) \} }(s) ds + \int_0^t b p |y_s|^{p} \mathbb{I}_{ \{ (0,r) \} }(s) dW_s \leq \\
|x_0|^p + \frac{b^2 p (p-1)}{2} \int_0^t |y_s|^{p}\mathbb{I}_{ \{
(0,r) \} }(s) ds+ \int_0^ t bp |y_s|^{p}\mathbb{I}_{ \{ (0,r) \}
}(s) dW_s
\end{eqnarray*}
Taking expectations we obtain,
\begin{eqnarray*}
\mathbb{E} |y_{t \wedge r}|^p \leq |x_0|^p + \frac{b^2 p (p-1)}{2}
\int_0^t \mathbb{E} |y_{s \wedge r}|^{p} ds.
\end{eqnarray*}

Therefore, using Gronwall inequality, we arrive at,
\begin{eqnarray}
\mathbb{E} |y_{t \wedge r}|^p \leq C(p) |x_0|^p ,
\end{eqnarray}
with $C(p)$ independent of $R > 0$. But $\mathbb{E} |y_{t \wedge
r}|^p = \mathbb{E} (|y_{t \wedge r}|^p \mathbb{I}_{ \{ r \geq t \}
}) +
 R^p P ( r < t  )$. That means that $P( t \wedge r < t) = P ( r < t  ) \to 0$ as
 $R \to \infty$ so $t \wedge r \to t$ in probability and noting that $r$ increases as $R$ increases we have that $t \wedge r \to t$
  almost surely too, as $R \to \infty$. Going back to (4) and
 using  Fatou's lemma  we obtain,
 \begin{eqnarray*}
\mathbb{E} |y_t|^p \leq C(p)|x_0|^p ,
\end{eqnarray*}
for any $p \geq 2$ and this is crucial to ensure that $\mathbb{E}
\int_0^t b p |y_{s}|^p dW_s = 0$ in the next step. Using Ito's
formula again on $|y_t|^p$, taking supremum and then expectations
and finally Doob's martingale inequality we arrive at,
\begin{eqnarray*}
\mathbb{E}  ( \sup_{0 \leq t \leq T} |y_t|^p) < A,
\end{eqnarray*}
for some $A > 0$.  To prove now that
$$\mathbb{E} ( \sup_{0 \leq t \leq T} |x_t|^p) < A$$ we use
Theorem 2.4.1 of \cite{Mao2} and obtain first the bound,
\begin{eqnarray*}
\mathbb{E}|x_t|^p \leq A,
\end{eqnarray*}
for any $p \geq 2$ and then using Ito's formula on $|x_t|^p$
taking supremum and then expectations and finally Doob's
martingale inequality we obtain the desired bound.

Finally, we shall prove that $$\mathbb{E}( \sup_{0 \leq t \leq T}
\frac{1}{x_t^2}) < A.$$ Set $l = \inf \{ t \geq 0 : \frac{1}{x_t}
\geq m \}$. Using Ito's formula on $(\frac{1}{x_{t \wedge l}})^2$ we
have,
\begin{eqnarray*}
\mathbb{E}(\frac{1}{x_{t \wedge l}})^2 = \mathbb{E}(\frac{1}{x_0})^2
+ \mathbb{E}\int_0^t (2  + 3 (\frac{b^2}{x_{s \wedge l}})^{2})ds.
\end{eqnarray*}
Using Gronwall's inequality we can prove that
$\mathbb{E}(\frac{1}{x_{t \wedge l}})^2 < A$ with $A$ independent of
$m$. As before we deduce that $\mathbb{E}(\frac{1}{x_{t}})^2 < A$
and then that $$\mathbb{E}( \sup_{0 \leq t \leq T} \frac{1}{x_t^2})
< A.$$

\qed

Now we can apply Theorem 1 to prove that our explicit numerical
scheme converges in the mean square sense. Let us note that from the
above moment bounds for the true solution we deduce that $x_t \geq
0$ a.s. Indeed, we have proved that $$\mathbb{E}(\frac{1}{x_{t
\wedge l}^2})  = \mathbb{E}(\frac{1}{x_{t \wedge l}^2} \mathbb{I}_{
\{l > t \} }) + m^2 P(l \leq t) < A.$$ Therefore, letting $m \to
\infty$ we obtain that $P( l \leq t) \to 0$, noting that $$P(x_t
\leq 0) = P(\bigcap_{m=1}^{\infty} \{ x_t \leq \frac{1}{m} \}) =
\lim_{n \to \infty} P(\{ x_t \leq \frac{1}{m} \}) \leq \lim_{m \to
\infty} P(l \leq t)=0.$$ Thus, we need our numerical scheme to be
positivity preserving. As we can see easily our scheme has this
advantage.
 For this kind of
problems there is also the  tamed Euler scheme as the authors in
\cite{Jentzen} proposes. This kind of method behaves very well and
in comparison with our ''semi-discrete'' method seems that the
tamed Euler method is less expensive because we have in each step
to compute an exponential. However, the tamed Euler scheme does
not seem to preserve positivity, at least,  for any $\Delta > 0$.

We shall give a simulation to show that the two numerical schemes
are close. So, let the following SDE,
\begin{eqnarray*}
x_t = 1 + \int_0^t -x_s^3 ds + \int_0^t x_s dW_s, \quad t \in
[0,1].
\end{eqnarray*}
We apply the tamed Euler scheme and the semi discrete scheme to
this problem for $\Delta = 10^{-4}$ and we plot the difference
between these methods, i.e. if $z_t = y^{tamed}_t-y^{semi}_t$ we
plot $z_t$ on $[0,1]$.

\begin{figure}[h!]
  \caption{  Difference between the semi discrete scheme and tamed Euler scheme for $x_0 = 1$, $\Delta=10^{-4}$, $b=1$, $T=1$.}
  \centering
    \includegraphics[width=0.5\textwidth]{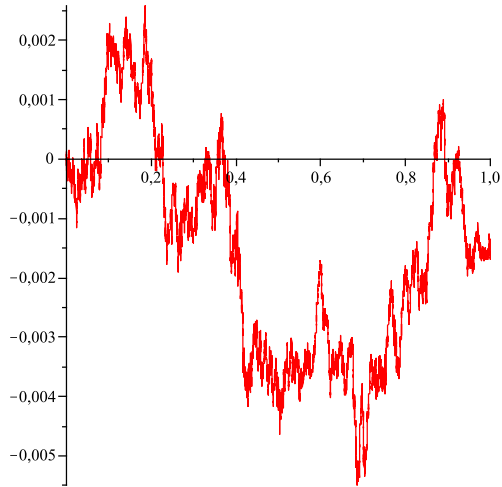}
\end{figure}

\begin{figure}[h!]
  \caption{ Tamed Euler method does not preserve positivity, $x_0 = 1$, $\Delta=10^{-3}$, $b=20$, $T=1$.}
  \centering
    \includegraphics[width=0.5\textwidth]{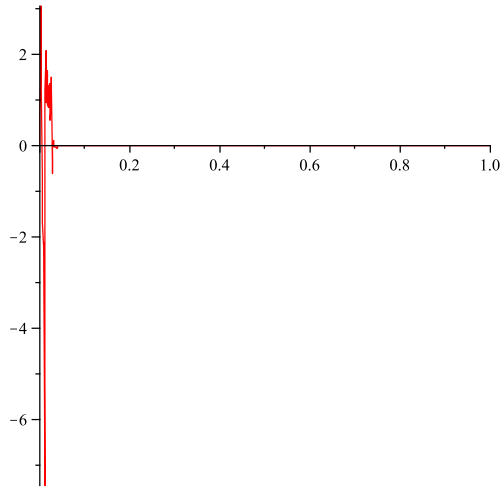}
\end{figure}
\section{Conclusion}

In this paper we propose a new numerical method for solving
stochastic differential equations. We apply our method to a
super-linear SDE and compare with the tamed Euler method. We see
that our method preserves positivity of the true solution. Our
method seems to behave very well when applied on super-linear
problems. In order to manage other interesting problems a suitable
transformation  and then the application of this ''semi-discrete''
method maybe the answer. Let us point out that our main result
remains true also for a multidimensional case.   Our goal in the
future is to apply our method to other stochastic differential
equations arising in financial mathematics and to give a more
detailed study of this method studying the rate of convergence, the
case where the choice $f,g$ are not locally Lipschitz etc.

\end{document}